\input amstex\documentstyle{amsppt}  
\pagewidth{12.5cm}\pageheight{19cm}\magnification\magstep1
\topmatter
\title An involution based left ideal in the Hecke algebra\endtitle
\author G. Lusztig\endauthor
\address{Department of Mathematics, M.I.T., Cambridge, MA 02139}\endaddress
\thanks{Supported in part by National Science Foundation grant 1303060.}\endthanks  \endtopmatter   
\document

\define\Irr{\text{\rm Irr}}

\define\mpb{\medpagebreak}

\define\bL{\bar L}

\define\hfH{\hat{\fH}}

\define\si{\sim}

\define\lb{\linebreak}

\define\op{\oplus}
   
\define\part{\partial}
\define\emp{\emptyset}
\define\imp{\implies}

\define\n{\notin}

\define\m{\mapsto}
\define\do{\dots}

\define\sub{\subset}    

\define\T{\times}
\define\ti{\tilde}
\define\nl{\newline}
\redefine\i{^{-1}}
\define\fra{\frac}
\define\un{\underline}
\define\ov{\overline}
\define\ot{\otimes}

\define\Hom{\text{\rm Hom}}
\define\End{\text{\rm End}}

\define\tr{\text{\rm tr}}

\define\a{\alpha}

\redefine\c{\chi}
\define\g{\gamma}
\redefine\d{\delta}
\define\e{\epsilon}

\define\p{\pi}
\define\ph{\phi}

\define\r{\rho}
\define\s{\sigma}
\redefine\t{\tau}

\define\k{\kappa}
\redefine\l{\lambda}

\define\x{\xi}

\redefine\G{\Gamma}

\define\CC{\bold C}

\define\FF{\bold F}

\define\II{\bold I}
\define\JJ{\bold J}

\define\NN{\bold N}

\define\QQ{\bold Q}

\define\ZZ{\bold Z}

\define\ca{\Cal A}
\define\cb{\Cal B}

\define\cd{\Cal D}
\define\ce{\Cal E}
\define\cf{\Cal F}

\define\ch{\Cal H}

\define\cm{\Cal M}

\define\cv{\Cal V}

\define\cx{\Cal X}

\define\fH{\frak H}

\define\fM{\frak M}

\define\ts{\ti s}

\define\tL{\ti L}

\define\Mod{\text{\rm Mod}}

\define\HZ{HZ}
\define\GE{Ge}
\define\KO{Ko}
\define\ORA{L1}
\define\LEA{L2}
\define\HEC{L3}
\define\INV{L4}
\define\BAR{L5}
\define\POWER{L6}
\define\LV{LV}

\head Introduction\endhead
\subhead 0.1\endsubhead
Let $W$ be a Coxeter group with set of simple reflections $S$ and with length function $l:W@>>>\NN$. 
Let $u$ be an indeterminate. Let $\fH$ be the $\QQ(u)$-vector space with basis $\{T_w;w\in W\}$. We regard 
$\fH$ as an associative $\QQ(u)$-algebra (with $1$) with multiplication defined by 
$T_wT_{w'}=T_{ww'}$ if $l(ww')=l(w)+l(w')$, $(T_s+1)(T_s-u^2)=0$ if $s\in S$. 
Let $*:W@>>>W$ (or $w\m w^*$) be an automorphism of $W$ such that $S^*=S$, $*^2=1$. Let 
$\II_*=\{w\in W;w^*=w\i\}$ be the set of "twisted involutions" of $W$.
 Let $M$ be the $\QQ(u)$-vector space with basis $\{a_w;w\in\II_*\}$. Following 
\cite{\LV}, for any $s\in S$ we define a $\QQ(u)$-linear map $T_s:M@>>>M$ by

$T_sa_w=ua_w+(u+1)a_{sw}$ if $sw=ws^*>w$;

$T_sa_w=(u^2-u-1)a_w+(u^2-u)a_{sw}$ if $sw=ws^*<w$;

$T_sa_w=a_{sws^*}$ if $sw\ne ws^*>w$;

$T_sa_w=(u^2-1)a_w+u^2a_{sws^*}$ if $sw\ne ws^*<w$.
\nl
(For $x,y$ in $W$ such that $y\i x\in S$ or $xy\i\in S$ we write $x<y$ or $y>x$ instead of $l(x)=l(y)-1$.)
According to \cite{\LV} and \cite{\BAR}, these linear maps define an $\fH$-module structure on $M$.
Let $\hfH$ be the vector space consisting of all formal (possibly infinite) sums 
$\sum_{x\in W}c_xT_x$ where $c_x\in\QQ(u)$. We can view $\fH$ as a subspace of $\hfH$ in an obvious way.
The $\fH$-module structure on $\fH$ (left multiplication) extends in an obvious way to an $\fH$-module 
structure on $\hfH$. We set
$$X=\sum_{x\in W;x^*=x}u^{-l(x)}T_x\in\hfH.$$
The following is the main result of this paper:

\proclaim{Theorem 0.2} (a) There exists a unique $\fH$-linear map $\mu:M@>\si>>\hfH$ such that 
$\mu(a_1)=X$. Moreover, $\mu$ is an isomorphism of $M$ onto the $\fH$-submodule of $\hfH$ generated by $X$.

(b) Let $z\in\II_*$; we set $\mu(a_z)=\sum_{x\in W}N^x_zT_x$ where $N^x_z\in\QQ(u)$. For any 
$x\in W$ we have $N^x_z\in\ZZ[u\i]$, hence we can define $n^x_z=N^x_z|_{u\i=0}\in\ZZ$.

(c) There is a unique surjective function $\p:W@>>>\II_*$ such that for $x\in W$, $z\in\II_*$ we have
$n^x_z=1$ if $z=\p(x)$, $n^x_z=0$ if $z\ne\p(x)$. (Note that $\p(1)=1$.)
\endproclaim
This was conjectured in \cite{\INV, 3.4, 3.7} where it was verified for several $W$ of low rank. In the case
where $W$ is of type $A_n$ and $*=1$, part (a) of the theorem was first proved by Hu and Zhang \cite{\HZ}. 
The proof of the theorem is given in Section 1.
In Section 2 we will discuss a special case of Theorem 0.2.
Section 3 is preparatory for Section 4.
In Section 4 we discuss some applications of Theorem 0.2. For example, we show that
if $W$ is a Weyl group of type $A_n$ and if $E$ is an irreducible representation of $\fH$,
then the action of $X$ on $E$ is through an operator of rank $1$; in particular the image of this operator
is a canonical line in $E$. As another application we show that if $W$ is a Weyl 
group of classical type and $E$ is an irreducible special representation of the 
asymptotic Hecke algebra attached to $W$ then
$E$ admits a basis such that any canonical basis element of that algebra acts in 
this basis through a matrix with all entries in $\NN$. 
A third application is a definition of a canonical $G(\FF_q)$-stable
subpace $\cf'$ of the space of 
functions $\cf$ on the flag manifold of a Chevalley group $G(\FF_q)$
over a finite field $\FF_q$ with the following properties: if $G=SL_n$, then
$\cf'$ contains exactly one copy of each irreducible representation of $G(\FF_q)$
which appears
in $\cf$; in general, the dimension of $\cf'$ is a polynomial in $q$ with
coefficients in $\NN$ whose value at $1$ is the number of involutions in $W$.
This polynomial is the sum of the fake degrees of the various irreducible
representations of the Hecke algebra (each one taken once). 

\head 1. Proof of Theorem 0.2\endhead
\subhead 1.1\endsubhead
The $\ZZ[u]$-submodule of $M$ with basis $\{a_w;w\in\II_*\}$ is stable under the maps $T_s:M@>>>M$
($s\in S$) hence is stable under the action of $T_x$ ($x\in W$) since $T_x$ is a composition of various $T_s$.
Hence for $x\in W$ we can write uniquely
$$T_xa_1=\sum_{z\in\II_*}L^x_za_z$$
where $L^x_z\in\ZZ[u]$. 

\subhead 1.2\endsubhead 
For $x\in W,z\in\II_*$, $s\in S$ we show:

$(u^2-u)L^x_{sz}=u^2L^{sx}_z+(u^2-u-1)L^x_z$ if $sz=zs^*>z, sx<x$;
    
$(u+1)L^x_{sz}-uL^x_z=u^2L^{sx}_z$ if $sz=zs^*<z, sx<x$;
    
$u^2L^x_{szs^*}=u^2L^{sx}_z+(u^2-1)L^x_z$ if $sz\ne zs^*>z, sx<x$;
  
$L^x_{szs^*}=u^2L^{sx}_z$ if $sz\ne zs^*<z,sx<x$;  

$uL^x_z+(u^2-u)L^x_{sz}=L^{sx}_z$ if $sz=zs^*>z,sx>x$;
    
$(u+1)L^x_{sz}+(u^2-u-1)L^x_z=L^{sx}_z$ if $sz=zs^*<z,sx>x$;
     
$u^2L^x_{szs^*}=L^{sx}_z$ if $sz\ne zs^*>z, sx>x$; 

$L^x_{szs^*}+(u^2-1)L^x_z=L^{sx}_z$ if $sz\ne zs^*<z,sx>x$.

\subhead 1.3\endsubhead 
For $x\in W,s\in S$ we have $T_sT_xa_1=\sum_{z\in\II_*}L^x_zT_sa_z$.
Note that $T_sT_xa_1=T_{sx}a_1$ if $sx>x$ and  $T_sT_xa_1=u^2T_{sx}a_1+(u^2-1)T_xa_1$ if $sx<x$. Thus,

$\sum_{z\in\II_*}L^x_zT_sa_z=\sum_{z\in\II_*}L^{sx}_za_z$ if $sx>x$,

$\sum_{z\in\II_*}L^x_zT_sa_z=\sum_{z\in\II_*}u^2L^{sx}_za_z+\sum_{z\in\II_*}(u^2-1)L^x_za_z$ if $sx<x$.  

Using the formulas for $T_sa_z$ in 0.1 we see that
$$\align&\sum_{z\in\II_*;sz=zs^*,sz>z}L^x_z(ua_z+(u+1)a_{sz})\\&+
\sum_{z\in\II_*;sz=zs^*,sz<z}L^x_z((u^2-u-1)a_z+(u^2-u)a_{sz})\\&
+\sum_{z\in\II_*;sz\ne zs^*,sz>z}L^x_za_{szs^*}
+\sum_{z\in\II_*;sz\ne zs^*,sz<z}L^x_z((u^2-1)a_z+u^2a_{szs^*})\endalign$$
or equivalently
$$\align&\sum_{z\in\II_*;sz=zs^*,sz>z}uL^x_za_z+\sum_{z\in\II_*;sz=zs^*,sz<z}(u+1)L^x_{sz}a_z\\&
+\sum_{z\in\II_*;sz=zs^*,sz<z}(u^2-u-1)L^x_za_z+\sum_{z\in\II_*;sz=zs^*,sz>z}(u^2-u)L^x_{sz}a_z\\&
+\sum_{z\in\II_*;sz\ne zs^*,sz<z}L^x_{szs^*}a_z+\sum_{z\in\II_*;sz\ne zs^*,sz<z}(u^2-1)L^x_za_z\\&
+\sum_{z\in\II_*;sz\ne zs^*,sz>z}u^2L^x_{szs^*}a_z\endalign$$
is equal to 

$\sum_{z\in\II_*}L^{sx}_za_z$ if $sx>x$, 

or to

$\sum_{z\in\II_*}u^2L^{sx}_za_z+\sum_{z\in\II_*}(u^2-1)L^x_za_z$ if $sx<x$.
\nl
We now take the coefficients of $a_z$ in the two sides of this equality. We obtain the equalities in 1.2.

\subhead 1.4\endsubhead 
Let $\bar{}:\ZZ[u,u\i]@>>>\ZZ[u,u\i]$ be the ring involution such that $\ov{u^n}=(-u)^{-n}$ for any $n\in\ZZ$.
We apply $\bar{}$ to the equalities in 1.2 and multiply the resulting equalities by $u^2$. We obtain the
following equalities.

$(1+u)\bL^x_{sz}=\bL^{sx}_z+(-u^2+u+1)\bL^x_z$ if $sz=zs^*>z, sx<x$;
    
$(u^2-u)\bL^x_{sz}+u\bL^x_z=\bL^{sx}_z$ if $sz=zs^*<z, sx<x$;
    
$\bL^x_{szs^*}=\bL^{sx}_z+(1-u^2)\bL^x_z$ if $sz\ne zs^*>z, sx<x$;
  
$u^2\bL^x_{szs^*}=\bL^{sx}_z$ if $sz\ne zs^*<z,sx<x$;  

$-u\bL^x_z+(1+u)\bL^x_{sz}=u^2\bL^{sx}_z$ if $sz=zs^*>z,sx>x$;    

$(u^2-u)\bL^x_{sz}+(-u^2+u+1)\bL^x_z=u^2\bL^{sx}_z$ if $sz=zs^*<z,sx>x$;
     
$\bL^x_{szs^*} =u^2\bL^{sx}_z$ if $sz\ne zs^*>z, sx>x$; 

$u^2\bL^x_{szs^*}+(1-u^2)\bL^x_z=u^2\bL^{sx}_z$ if $sz\ne zs^*<z,sx>x$.  

\subhead 1.5\endsubhead 
It is well known (see for example \cite{\BAR}) that there is a unique function $\ph:\II_*@>>>\NN$ such that
$\ph(1)=0$ and such that for any $z\in\II_*$ and any $s\in S$ such that $sz<z$ we have
$\ph(z)=\ph(sz)+1$ if $sz=zs^*$ and $\ph(z)=\ph(szs^*)$ if $sz\ne zs^*$.
By induction on $l(z)$ we see that $\ph(z)=l(z)\mod2$ for any $z\in\II_*$. Hence for $z\in\II_*$ we can set 
$\e(z)=(-1)^{(l(z)+\ph(z))/2}$. From the definitions we see that for any $z\in\II_*$ and any $s\in S$ we have

(a) $\e(z)=-\e(sz)$ if $sz=zs^*$ and $\e(z)=-\e(szs^*)$ if $sz\ne zs^*$.

\subhead 1.6\endsubhead 
For $x\in W,z\in\II_*$ we set
$$\tL^x_z=(-1)^{l(x)}\e(z)\bL^x_z.$$ 
With this notation the formulas in 1.4 can be rewitten as follows.

$(1+u)\tL^x_{sz}=\tL^{sx}_z+(u^2-u-1)\tL^x_z$ if $sz=zs^*>z, sx<x$;
    
$(u^2-u)\tL^x_{sz}-u\tL^x_z=\tL^{sx}_z$ if $sz=zs^*<z, sx<x$;
    
$\tL^x_{szs^*}=\tL^{sx}_z+(u^2-1)\tL^x_z$ if $sz\ne zs^*>z, sx<x$;
  
$u^2\tL^x_{szs^*}=\tL^{sx}_z$ if $sz\ne zs^*<z,sx<x$;  

$u\tL^x_z+(1+u)\tL^x_{sz}=u^2\tL^{sx}_z$ if $sz=zs^*>z,sx>x$;    

$(u^2-u)\tL^x_{sz}+(u^2-u-1)\tL^x_z=u^2\tL^{sx}_z$ if $sz=zs^*<z,sx>x$;
     
$\tL^x_{szs^*} =u^2\tL^{sx}_z$ if $sz\ne zs^*>z, sx>x$; 

$u^2\tL^x_{szs^*}+(u^2-1)\tL^x_z=u^2\tL^{sx}_z$ if $sz\ne zs^*<z,sx>x$.  

\subhead 1.7\endsubhead 
Giving an $\fH$-linear map $\mu:M@>>>\hfH$ is the same as giving a family of elements
$Y_z\in\hfH$ (one for each $z\in\II_*$) such that for any $z\in\II_*,s\in S$ we have

$T_sY_z=uY_z+(u+1)Y_{sz}$ if $sz=zs^*>z$;

$T_sY_z=(u^2-u-1)Y_z+(u^2-u)Y_{sz}$ if $sz=zs^*<z$;

$T_sY_z=Y_{szs^*}$ if $sz\ne zs^*>z$;

$T_sY_z=(u^2-1)Y_z+u^2Y_{szs^*}$ if $sz\ne zs^*<z$.

Indeed, if $\mu$ is given then the elements $Y_z=\mu(a_z)$ satisfy the equations above. Conversely, if
the elements $Y_z$ are given as above then we can define a $\QQ(u)$-linear map $\mu:M@>>>\hfH$ by
$\mu(a_z)=Y_z$ for all $z\in\II_*$. 
This map will be compatible with the action of $T_s$ for any $s\in S$ hence it will be
automatically $\fH$-linear. Setting $Y_z=\sum_{x\in W}N^x_zT_z$ where $N^x_z\in\QQ(u)$ we see that 
giving an $\fH$-linear map $\mu:M@>>>\hfH$ is the same as giving a family of elements
$\{N^x_z;(x,z)\in W\T\II_*\}$ in $\QQ(u)$ such that the following equations are satisfied for any 
$z\in\II_*$, $s\in S$:

$\sum_{x\in W}N^x_zT_sT_x=\sum_{x\in W}uN^x_zT_x+\sum_{x\in W}(u+1)N^x_{sz}T_x$ if $sz=zs^*>z$;

$\sum_{x\in W}N^x_zT_sT_x=\sum_{x\in W}(u^2-u-1)N^x_zT_x+\sum_{x\in W}(u^2-u)N^x_{sz}T_x$ if $sz=zs^*<z$;

$\sum_{x\in W}N^x_zT_sT_x=\sum_{x\in W}N^x_{szs^*}T_x$ if $sz\ne zs^*>z$;

$\sum_{x\in W}N^x_zT_sT_x=\sum_{x\in W}(u^2-1)N^x_zT_x+\sum_{x\in W}u^2N^x_{szs^*}T_x$ if $sz\ne zs^*<z$.
\nl
(We then say that the family $\{N^x_z;(x,z)\in W\T\II_*\}$ is admissible.) Here we replace
$$\align&\sum_{x\in W}N^x_zT_sT_x=\sum_{x\in W;sx>x}N^x_zT_{sx}+\sum_{x\in W;sx<x}u^2N^x_zT_{sx}+
\sum_{x\in W;sx<x}(u^2-1)N^x_zT_x\\&=
\sum_{x\in W;sx<x}N^{sx}_zT_x+\sum_{x\in W;sx>x}u^2N^{sx}_zT_x+\sum_{x\in W;sx<x}(u^2-1)N^x_zT_x\\&
=\sum_{x\in W;sx<x}(N^{sx}_z+(u^2-1)N^x_z)T_x+\sum_{x\in W;sx>x}u^2N^{sx}_zT_x.\endalign$$
We see that the condition that $\{N^x_z;(x,z)\in W\T\II_*\}$ is admissible is equivalent to the following
set of equations (with $x\in W,z\in\II_*,s\in S$).

$(1+u)N^x_{sz}=N^{sx}_z+(u^2-u-1)N^x_z$ if $sz=zs^*>z, sx<x$;
    
$(u^2-u)N^x_{sz}-uN^x_z=N^{sx}_z$ if $sz=zs^*<z, sx<x$;
    
$N^x_{szs^*}=N^{sx}_z+(u^2-1)N^x_z$ if $sz\ne zs^*>z, sx<x$;
  
$u^2N^x_{szs^*}=N^{sx}_z$ if $sz\ne zs^*<z,sx<x$;  

$uN^x_z+(1+u)N^x_{sz}=u^2N^{sx}_z$ if $sz=zs^*>z,sx>x$;    

$(u^2-u)N^x_{sz}+(u^2-u-1)N^x_z=u^2N^{sx}_z$ if $sz=zs^*<z,sx>x$;
     
$N^x_{szs^*} =u^2N^{sx}_z$ if $sz\ne zs^*>z, sx>x$; 

$u^2N^x_{szs^*}+(u^2-1)N^x_z=u^2N^{sx}_z$ if $sz\ne zs^*<z,sx>x$.  
\nl
Comparing with the formulas in 1.6, we see that the family $\{\tL^x_z;(x,z)\in W\T\II_*\}$ is admissible. 
Hence there is a unique $\fH$-linear map $\mu:M@>>>\hfH$ such that for any $z\in\II_*$ we have 
$\mu(a_z)=\sum_{x\in W}\tL^x_zT_z$. Since $L^x_1=\d_{x,x^*}u^{l(x)}$ (see \cite{\POWER, 1.8}), we have
$\tL^x_1=\d_{x,x^*}(-1)^{l(x)}(-u)^{-l(x)}=\d_{x,x^*}u^{-l(x)}$, so that $\mu(a_1)=X$ (see 0.1). Thus the
existence part in 0.2(a) is established. The uniqueness part in 0.2(a) is obvious since $a_1$ generates $M$
as an $\fH$-module. Since $L^x_z\in\ZZ[u]$ we see that $\tL^x_z\in\ZZ[u\i]$ and 0.2(b) is established.

\subhead 1.8\endsubhead 
The algebra $\fH$ and its module $M$ can be specialized to $u=0$. Then
$\fH$ becomes a $\QQ$-algebra $\fH_0$ with basis $\{\un T_w;w\in W\}$ and multiplication given by
$\un T_w\un T_{w'}=\un T_{ww'}$ if $l(ww')=l(w)+l(w')$, $(\un T_s+1)\un T_s=0$ if $s\in S$; $M$ becomes
a $\QQ$-vector space $M_0$ with basis $\{\un a_w;w\in\II_*\}$ and with $\fH_0$-module structure given by

$\un T_s\un a_w=\un a_{sw}$ if $sw=ws^*>w$;

$\un T_s\un a_w=\un a_{sws^*}$ if $sw\ne ws^*>w$;

$\un T_s\un a_w=-a_w$ if $sw<w$.
\nl
Here $s\in S,w\in\II_*$. We have the following result:

(a) {\it There is a unique map $W\T\II_*@>>>\II_*$, $(x,w)\m x\circ w$ such that
$\un T_x\un a_w=\e_{x,w}\un a_{x\circ w}$ for any $x\in W,w\in\II_*$; here $\e_{x,w}=\pm1$ is a well
defined sign.}
\nl
We argue by induction on $l(x)$. If $x=1$ we have
$\un T_x\un a_w=\un a_w$ so that $1\circ w=w,\e_{1,w}=1$. Assume now that $x\ne1$. Let $s\in S$, $x'\in W$
 be such that
$x=sx'>x'$. By the induction hypothesis we have $\un T_{x'}\un a_w=\pm\un a_{u'}$ for some $u'\in\II_*$.
Hence $\un T_x\un a_w=\pm \un T_s\un a_{u'}$ and this equals $\pm\un a_u$ for some $u\in\II_*$. This proves
(a).

We show:

(b) $\un T_x\un a_1=(-1)^{l(x)}\e(x\circ 1)\un a_{x\circ 1}$ for any $x\in W$.
\nl
We argue by induction on $l(x)$. If $x=1$ we have 
$\un T_x\un a_1=\un a_1$ hence $1\circ1=1$ and $l(x)=0$, $\e(1\circ1)=\e(1)=1$ and the result holds.
Assume now that $x\ne1$. Let $s\in S$, $x'\in W$
 be such that $x=sx'>x'$. By the induction hypothesis we have 
$\un T_{x'}\un a_1=(-1)^{l(x')}\e(w)\un a_w$ where $w=x'\circ1$.
We have $\un T_x\un a_1=\un T_s\un T_{x'}\un a_1=(-1)^{l(x')}\e(w)\un T_s\un a_w$.
Now $\un T_s\un a_w=f(s,w)\un a_{s\circ w}$ where $f(s,w)=1$ if $sw>w$, $f(s,w)=-1$ if $sw<w$. 
It is enough to prove that $(-1)^{l(x)}\e(x\circ1)=(-1)^{l(x')}\e(w)f(s,w)$. Since $l(x)=l(x')+1$ it is
enough to prove that $\e(x\circ1)=-\e(w)f(s,w)$. We have $x\circ1=s\circ w$ hence it is enough to prove that
$\e(s\circ w)=-\e(w)f(s,w)$ or that:
$\e(sw)=-\e(w)$ if $sw=ws^*>w$, $\e(sws^*)=-\e(w)$ if $sw\ne ws^*>w$.
This is clear from the definition of $\e$. This proves (b).

\mpb

We define $\p:W@>>>\II_*$ by $\p(x)=x\circ 1$. We show:

(c) {\it $\p$ is surjective.}
\nl
Let $w\in\II_*$. We show by induction on $l(w)$ that $w\in\p(W)$. If $w=1$ we have $w=\p(1)$. Assume now
that $w\ne1$. We can find $s\in S$ such that $sw<w$.
Assume first that $sw=ws^*$. Then by the induction hypotesis we have $sw=x\circ 1$ for some $x\in W$ hence
$\un T_s\un T_x\un a_1=\pm\un T_s\un a_{sw}=\pm\un a_w$; moreover $\un T_s\un T_x$ equals $\un T_x$
(if $sx>x$) or $\un T_{sx}$ (if $sx<x$). Thus $w=x\circ 1$ or $w=sx\circ 1$.

Assume next that $sw\ne ws^*$. Then by the induction hypotesis we have $sws^*=x\circ 1$ for some $x\in W$ 
hence $\un T_s\un T_x\un a_1=\pm\un T_s\un a_{sws^*}=\pm\un a_w$; moreover $\un T_s\un T_x$ equals $\un T_x$
(if $sx>x$) or $\un T_{sx}$ (if $sx<x$). Thus $w=x\circ 1$ or $w=sx\circ 1$. This proves (c).

\subhead 1.9\endsubhead 
For $x\in W$ we have $\un T_x\un a_1=\sum_{z\in \II_*}\un L^x_z\un a_z$
where $\un L^x_z=L^x_z|_{u=0}\in\ZZ$. Comparing with 1.8(b) we see that
$\un L^x_z=(-1)^{l(x)}\e(z)$ if $z=x\circ 1$  and $\un L^x_z=0$ if $z\ne x\circ 1$.
It follows that $\tL^x_z|_{u\i=0}=1$ if $z=x\circ 1$ and $\tL^x_z|_{u\i=0}=0$ if $z\ne x\circ 1$.
Thus 0.2(c) holds.

\subhead 1.10\endsubhead 
We show that the map $\mu:M@>>>\hfH$ is injective. It is enough to show that
the elements $\{\mu(a_z);z\in\II_*\}$ are linearly independent.
Assume that $\sum_{z\in\II_*}\x_z\mu(a_z)=0$ where $\x_z\in\QQ(u)$ are zero for all but finitely many $z$
and $\x_z\ne0$ for some $z\in\II_*$. We can assume that $\x_z\in\ZZ[u\i]$ for all $z$ and 
$\x_z|_{u\i=0}\ne0$ for some $z=z_0$. We have $\sum_{z\in\II_*,x\in W}\x_z\tL^x_zT_x=0$ hence 
$\sum_{z\in\II_*}\x_z\tL^x_z=0$ for any $x\in W$. Setting $u\i=0$ we deduce that
$\sum_{z\in\II_*}\x_z|_{u\i=0}n^x_z=0$ for any $x\in W$. By 0.2(c) this can be written as
$\x_{\p(x)}|_{u\i=0}=0$ for any $x\in W$. By 1.8(c) we can find $x\in W$ such that
$\p(x)=z_0$. For this $x$ we have $\x_{z_0}|_{u\i=0}=0$. This is a contradiction, Thus the injectivity
of $\mu$ is proved. This completes the proof of Theorem 0.2.

\subhead 1.11\endsubhead
In the case where $W$ is of type $A_1$ with $S=\{s\}$ we have
$\mu(a_1)=u\i T_s+1$, $\mu(a_s)=(u-1)u\i T_s$. 

In the case where $W$ is of type $A_2$ with $S=\{s,t\}$ and $*=1$ we have

$\mu(a_1)=u^{-3}T_{sts}+u^{-2}T_{st}+u^{-2}T_{ts}+u\i T_s+u\i T_t+1$, 

$\mu(a_s)=(u-1)(u^{-3}T_{sts}+u^{-2}T_{st}+u\i T_s)$, 

$\mu(a_t)=(u-1)(u^{-3}T_{sts}+u^{-2}T_{ts}+u\i T_t)$, 

$\mu(a_{sts})=(u-1)((u\i+u^{-2}-u^{-3})T_{sts}+u\i T_{st}+u\i T_{ts})$.
\nl
(See \cite{\INV, 32, 3.3}.)

\subhead 1.12\endsubhead
For $x\in W,z\in\II_*$ we set $\tL^x_z=(u-1)^{\ph(z)}\l^x_z$ where $\ph(z)$ is as in 1.5 and 
$\l^x_z\in\QQ(u)$. We show:

(a) $\l^x_z\in\ZZ[u\i]$ and $\ov{\l^x_z}=(-u^2)^{l(x)+(1/2)(l(z)-\ph(z))}\l^x_z$.
\nl
From the definitions we have $L^1_z=\d_{1,z}$ hence $\tL^1_z=\d_{1,z}$ and $\l^1_z=\d_{1,z}$. From the
formulas in 1.6 (with $s\in S$) we deduce (assuming $sx<x$):

$\l^x_z=\l^x_{sz}-u\i \l^{sx}_z$ if $sz=zs^*<z$; 

$\l^x_z=u\i \l^{sx}_z+(1-u^{-2})\l^{sx}_{sz}$ if $sz=zs^*>z$; 

$\l^x_z=u^{-2}\l^{sx}_{szs^*}$ if $sz\ne zs^*>z$; 

$\l^x_z=\l^{sx}_{szs^*}+(1-u^{-2})\l^{sx}_z$ if $sz\ne zs^*<z$. 

From this (a) follows by induction on $l(x)$.

\subhead 1.13\endsubhead
In this subsection we give an application of the function $\e:\II_*@>>>\{\pm1\}$ in 1.5.
Let $E=\QQ(u)$ viewed as an $\fH$-module in which $T_x$ ($x\in W$) acts as multiplication by
$(-1)^{l(x)}$ (sign representation of $\fH$).
 We define a $\QQ(u)$-linear map $f:M@>>>E$ by $f(a_z)=\e(z)$. We claim that $f$ is 
$\fH$-linear. It is enough to show that for any $w\in\II_*,s\in S$ we have:

$-\e(w)=u\e(w)+(u+1)\e(sw)$ if $sw=ws^*>w$;

$-\e(w)=(u^2-u-1)\e(w)+(u^2-u)\e(sw)$ if $sw=ws^*<w$;

$-\e(w)=\e(sws^*)$ if $sw\ne ws^*>w$;

$-\e(w)=(u^2-1)\e(w)+u^2\e(sws^*)$ if $sw\ne ws^*<w$.
\nl
This follows from 1.5(a).

\head 2. The biregular representation of $\fH$\endhead
\subhead 2.1\endsubhead
In this section we discuss the special case of Theorem 0.2 in the case where $W$ in 0.1 is replaced by 
$W^2=W\T W$, $S$ is replaced by $S^2=S\T\{1\}\cup\{1\}\T S$ and $*:W@>>>W$ is replaced by $*:W^2@>>>W^2$, 
$(x,y)\m(y,x)$. In this case we have $\II_*=\{(x,y)\in W^2;xy=1\}$.

We use the notation $\fH,\hfH$ in reference to $W$. Let $\fH^{(2)},\hfH^{(2)}$ be the objects analogous to 
$\fH,\hfH$ defined in terms of $W^2$ instead of $W$. Thus $\fH^{(2)}=\fH\ot\fH$ with basis
$\{T_w\ot T_{w'};(w,w')\in W^2\}$ which is the analogue for $\fH^{(2)}$ of the basis $\{T_w;w\in W\}$ of 
$\fH$. We write $T_w$ (resp. $T'_w$) instead of $T_w\ot 1$ (resp. $1\ot T_w$). Then the basis element 
$T_w\ot T_{w'}$ is actually the product in $\fH^{(2)}$ of $T_w$ and $T'_{w'}$ in either order. In our case 
we have $M=\fH$ with the basis $\{a_{w,w\i}=T_w;w\in W\}$ viewed as an $\fH^{(2)}$-module in which the 
action of $T_x\ot T_y$ is $T_r\m T_xT_rT_{y\i}$. (We refer to this as the {\it biregular representation}.)
We can view $\fH^{(2)}$ as a subspace of $\hfH^{(2)}$ in an obvious way.
The $\fH^{(2)}$-module structure on $\fH^{(2)}$ (left multiplication) extends in an obvious way to a 
$\fH^{(2)}$-module structure on $\hfH^{(2)}$. 
The following is a restatement of Theorem 0.2 in our case.

\proclaim{Corollary 2.2} (a) There exists a unique $\fH^{(2)}$-linear map $\mu:\fH@>>>\hfH^{(2)}$
(where $\fH,\hfH^{(2)}$ are viewed as $\fH^{(2)}$-modules as above) such that
$\mu(1)=\sum_{w\in W}u^{-2l(w)}T_w\ot T_w\in\hfH^{(2)}$. Moreover, $\mu$ is an isomorphism of $\fH$ onto
the $\fH^{(2)}$-submodule of $\hfH^{(2)}$ generated by $\mu(1)$.

(b) Let $z\in W$; we set $\mu(T_z)=\sum_{(x,y)\in W^2}N^{x,y}_zT_x\ot T_y$ where $N^{x,y}_z\in\QQ(u)$. For 
any $(x,y)\in W^2$ we have $N^{x,y}_z\in\ZZ[u\i]$, hence we can define $n^{x,y}_z=N^{x,y}_z|_{u\i=0}\in\ZZ$.

(c) There is a unique surjective function $\p:W^2@>>>W$ such that for $(x,y)\in W$, $z\in W$ we have
$n^{x,y}_z=1$ if $z=\p(x,y)$, $n^{x,y}_z=0$ if $z\ne\p(x,y)$. (Note that $\p(1,1)=1$.)
\endproclaim

In the remainder of this section we shall indicate a proof of a part of Corollary which is somewhat 
different from that of Theorem 0.2.

\subhead 2.3\endsubhead
Let $\t:\fH@>>>\QQ(u)$ be the $\QQ(u)$-linear map such that $\t(T_x)=0$ if $x\ne1$, $\t(T_1)=1$. For
$x,y\in W$ we have $\t(T_xT_y)=0$ if $xy\ne1$, $\t(T_xT_y)=u^{2l(x)}$ if $xy=1$. (See \cite{\HEC, 10.4(a)}.) 
It follows that for $x,y$ in $W$ we have
$$T_xT_y=\sum_{z\in W}\t(T_xT_yT_z)u^{-2l(z)}T_{z\i}.$$
Since $T_x\m T_{x\i}$ defines an antiautomorphism of $\fH$ we have 
$$T_{y\i,x\i}=\sum_{z\in W}\t(T_xT_yT_z)u^{-2l(z)}T_z=\sum_{z\in W}\t(T_{y\i}T_{x\i}T_{z\i})u^{-2l(z)}T_z$$
hence 
$$\t(T_xT_yT_z)=\t(T_{y\i}T_{x\i}T_{z\i})\tag a$$
for any $x,y,z$ in $W$. By \cite{\HEC, 10.4(b)} we have $\t(hh')=\t(h'h)$ for any $h,h'$ in $\fH$. In 
particular for $x,y,z$ in $W$ we have 
$$\t(T_xT_yT_z)=\t(T_yT_zT_x)=\t(T_zT_xT_y).\tag b$$

\proclaim{Lemma 2.4} (a) For $x,y,z$ in $W$ we set $p_{x,y,z}=\t(T_xT_yT_z)u^{-2l(z)-2l(y)}$,
$p'_{x,y,z}=\t(T_xT_yT_z)u^{-2l(x}$. We have $p_{x,y,z}\in d_{x,y,z}+u\i\ZZ[u\i]$, 
$p'_{x,y,z}\in d'_{x,y,z}+u\ZZ[u]$ where $d_{x,y,z}\in\{0,1\}$, $d'_{x,y,z}\in\{0,\pm1\}$.
Moreover $p'_{x,y,z}=(-1)^{l(x)+l(y)+l(z)}\ov{p_{x,y,z}}$ hence $d'_{x,y,z}=(-1)^{l(x)+l(y)+l(z)}d_{x,y,z}$.

(b) Let $y,z$ be in $W$. There is exactly one $x\in W$ (denoted by $y*z$) such that $d_{x,y,z}=1$ (or 
equivalently such that $d'_{x,y,z}=\pm1$). For all other $x$ we have $d_{x,y,z}=d'_{x,y,z}=0$.
\endproclaim
We argue by induction on $l(z)$. If $z=1$ we have 
$p_{x,y,z}=\t(T_xT_y)u^{-2l(y)}=\d_{x,y}$,
$p'_{x,y,z}=\t(T_xT_y)u^{-2l(x)}=\d_{x,y}$. Hence (a),(b) hold with $y*z=y\i$.
(We have $\d_{x,y}= (-1)^{l(x)+l(y)}\d_{x,y}$.)

Assume now that $l(z)\ge1$. We write $z=sz'$, $s\in S$, $l(z')=l(z)-1$. If $ys>y$ we have 
by the induction hypothesis
$$p_{x,y,z}=\t(T_xT_yT_sT_{z'})u^{-2l(z')-2l(y)-2}=p_{x,ys,z'}\in d_{x,ys,z'}+u\i\ZZ[u\i],$$
$$p'_{x,y,z}=\t(T_xT_yT_sT_{z'})u^{-2l(x)}=p'_{x,ys,z'}\in d'_{x,ys,z'}+u\ZZ[u]$$
hence the result holds: we have $d_{x,y,z}=d_{x,ys,z'}$, $d'{x,y,z}=d'_{x,ys,z'}$, $y*z=(ys)*(sz)$.
(We use that $(-1)^{l(x)+l(y)+l(z)}=(-1)^{l(x)+l(ys)+l(z')}$.)

If $ys<y$ we have by the induction hypothesis
$$\align&p_{x,y,z}=\t(T_xT_yT_sT_{z'})u^{-2l(z')-2l(y)-2}\\&=
p_{x,ys,z'}u^{-2}+p_{x,y,z'}u^{-2}(u^2-1)\in d_{x,y,z'}+u\i\ZZ[u\i],\endalign$$
$$\align&p'_{x,y,z}=\t(T_xT_yT_sT_{z'})u^{-2l(x)}\\&=p'_{x,ys,z'}u^2+p'_{x,y,z'}(u^2-1)
\in -d'_{x,y,z'}+u\ZZ[u]\endalign$$
hence the result holds: we have $d_{x,y,z}=d_{x,y,z'}$, $d'_{x,y,z}=-d'_{x,y,z'}$, $y*z=y*(sz)$.

(We use that 

$(-1)^{l(x)+l(y)+l(z)}=(-1)^{l(x)+l(ys)+l(z')}$, $(-1)^{l(x)+l(y)+l(z)}=-(-1)^{l(x)+l(y)+l(z')}$,

$\ov{u^2-1}=-(1-u^{-2})$.)

\subhead 2.5\endsubhead
For $a\in W$ we show
$$T_aX=T'_{a\i}X.\tag a$$
We have
$$T_aX=\sum_{w,z\in W}u^{-2l(w)-2l(z)}\t(T_aT_wT_z\i)T_zT'_w\in\hfH^{(2)},$$
$$T'_{a\i}X=\sum_{w,z\in W}u^{-2l(w)-2l(z)}\t(T_{a\i}T_wT_{z\i})T_wT'_z\in\hfH^{(2)}.$$
Making the change of variable $(w,z)\m(z,w)$ in the last sum we obtain
$$T'_{a\i}X=\sum_{w,z\in W}u^{-2l(w)-2l(z)}\t(T_{a\i}T_zT_{w\i})T_zT'_w.$$
It remains to show:
$$\t(T_{a\i}T_zT_{w\i})=\t(T_aT_wT_{z\i}).$$
Indeed, by 2.3(a) the left hand side is equal to $\t(T_{z\i}T_aT_w)$ and by 2.3(b) this is equal to the right
hand side.

\subhead 2.6\endsubhead
We give an alternative proof of the existence of $\mu$ in Corollary 2.2. For any $a\in W$ we set 
$X_a=T_aX=T'_{a\i}X\in\hfH^{(2)}$, see 2.5(a). Thus, $X_1=X$.
We define a $\QQ(u)$-linear map $\mu:\fH@>>>\hfH^{(2)}$ by $T_a\m X_a$ for all $a\in W$.
For $h\in\fH$, $r\in W$ we have $\mu(T_rh)=T_r\mu(h)$ (using the description $X_a=T'_{a\i}X$)
and $\mu(hT_{r\i})=T'_r\mu(h)$ (using the description $X_a=T_aX$). It follows that $\mu$ is
$\fH^{(2)}$-linear.

In our case $\p:W^2@>>>W$ is given by $\p(x,y)=(y*(x\i),(y*(x\i))\i$.

\subhead 2.7\endsubhead
In the case where $W$ is of type $A_1$ with $S=\{s\}$ we have

$\mu(T_1)=T_1\ot T_1+u^{-2}T_s\ot T_s$,

$\mu(T_s)=T_1\ot T_s+T_s\ot T_1+(1-u^{-2})T_s\ot T_s$.

\head 3. $\G$-equivariant vector bundles on $\G$\endhead
\subhead 3.1\endsubhead
Let $\G$ be a finite group. Let $K_\G(\G)$ be the Grothendieck group of
$\G$-equivariant (complex) vector bundles on $\G$ where $\G$ acts on $\G$ by 
conjugation. For $x\in\G$ let $\G_x=Z_\G(x)$ and let $\Irr\G_x$ be a set of
representatives for the isomorphism classes of irreducible representations of
$\G_x$ over $\CC$. 
For any $x\in\G$ and any $\r\in\Irr\G_x$ there is a unique
(up to isomorphism) $\G$-equivariant vector bundle $E_{x,\r}$ on $\G$ such that
the support of $E_{x,\r}$ is the conjugacy class of $x$ and is such 
that the action of $\G_x$ on the fibre of $E_{x,\r}$ is isomorphic to $\r$.
Let $\un\G$ be a set of representatives for the conjugacy classes in $\G$.
Let $\fM(\G)=\{(x,\r);x\in\un\G,\r\in\Irr \G_x\}$. The classes of $E_{x,\r}$ 
(with $(x,\r)\in\fM(\G)$) form a $\ZZ$-basis of $K_\G(\G)$.

Following Kottwitz \cite{\KO} we consider the element $\k\in K_\G(\G)$ defined by
$$\k=\sum_{(x,\r)\in\fM(\G)}\sum_{s\in\G;s^2=x}\fra{|Z_{\G_x}|}{|\G_x|}
(1:\r|_{Z_{\G_x}(s)}) E_{x,\r}$$
where $(1:\r|_{Z_{\G_x}(s)})$ denotes the multiplicity of the unit representation
of $Z_{\G_x}(s)$.

\proclaim{Proposition 3.2}Define $V=\a_!\CC$ where $\a:\G@>>>\G$ is $g\m g^2$. 
Note that $V$ is a $\G$-equivariant vector bundle on $\G$. We have $V=\k$ in 
$K_\G(\G)$.
\endproclaim
Let $\G^{(2)}=\{(g,h)\in\G\T\G;gh=hg\}$. For any $\G$-equivariant vector bundle 
$\cv$ on $\G$ we define $\ph_{\cv}:\G^{(2)}@>>>\CC$ as follows: 
$\ph_{\cv}(g,h)$ is trace of the action of $h$ on the fibre of $\cv$ at $g$.
For example, if $(x,\r)\in\fM(\G)$, we have
$$\ph_{E_{x,\r}}(g,h)=|\G_x|\i\sum_{a\in\G;aga\i=x}\tr(aha\i,\r).$$
Note that $\cv@>>>\ph_{\cv}$ induces an injective linear map from the vector space 
$\CC\ot K_\G(\G)$ into the vector spaces of functions $\G^{(2)}@>>>\CC$, see 
\cite{\LEA}. Hence it suffices to show that $\ph_V=\ph_\k$. For $(g,h)\in\G^{(2)}$
we have
$$\align&
\ph_\k(g,h)=\sum_{x\in\G,\r\in\Irr\G_x}\fra{|\G_x|}{|\G|}\sum_{s\in\G;s^2=\x}
\fra{|Z_{\G_x}(s)|}{|\G_x|}(1:\r|_{\G_x\cap\G_s})\ph_{E_{x,\r}}(g,h)\\&
=\sum_{x\in\G,\r\in\Irr\G_x}|\G|\i\sum_{s\in\G;s^2=x}
\sum_{u\in\G_x\cap\G_s}\tr(u\i,\r)|\G_x|\i\sum_{a\in\G; aga\i=x} \tr(aha\i,\r)\\&
=\sum_{x\in\G}|\G|\i\sum_{s\in\G;s^2=x}
\sum_{u\in\G_x\cap\G_s}|\G_x|\i\sum_{a\in\G;aga\i=x}|\{z\in G_x;zaha\i z\i=u\i\}|
\endalign$$
Setting $s'=a\i sa, u'=a\i ua,z'=a\i za$ we obtain
$$\align&\ph_\k(g,h)=\sum_{s'\in\G;s'{}^2=g}
\sum_{u'\in\G_g\cap\G_{s'}}
|\G_g|\i|\{z'\in G_g;z'hz'{}\i=u'{}\i\}|\\&
=\sum_{s'\in\G;s'{}^2=g}|\G_g|\i|\{z'\in G_g;z'hz'{}\i\in G_g\cap G_{s'}\}|.\endalign$$
Sertting $\ts=z'{}\i s'z'$ we obtain
$$\ph_\k(g,h)=|\{\ts\in\G;\ts^2=g,\ts h=h\ts\}|.$$
From the definitions we have
$$\ph_V(g,h)=|\{\ts\in\G;\ts^2=g,\ts h=h\ts\}|.$$
The proposition is proved.

\subhead 3.3\endsubhead
As in \cite{\LEA, 2.5}, any $(y,\s)\in\fM(\G)$ defines a $\CC$-linear function
\lb $\c_{y,\s}:\CC\ot K_\G(\G)@>>>\CC$ by the rule
$$\c_{y,\s}(U)=(\dim\s)\i\sum_{\g\in G_y}\tr(y,U_\g)\tr(\g,\s)$$
for any $\G$-equivariant vector bundle. (This is in fact an algebra homomorphism 
for the algebra structure defined in \cite{\LEA, 2.2}.)
Moreover, if $(x,\r)\in M(\G)$ then
$$\c_{y,\s}(E_{x,\r})=\fra{|\G_y|}{\dim\s}\{(x,\r),(y,\s^*)\}\tag a$$
where $\{,\}$ is the nonabelian Fourier transform matrix of \cite{\ORA} and
$\s^*\in\Irr\G_y$ is isomorphic to the dual of $\s$.
We compute $\c_{y,\s}(V)$ where $V$ is as in 3.2 and $\s$ has
Frobenius-Schur indicator $1$. By the proof of 3.2 we have
$$\align&\c_{y,\s}(V)=(\dim\s)\i\sum_{\g\in\G_y}\tr(y,V_\g)\tr(\g,\s)\\&
=(\dim\s)\i\sum_{\g\in\G_y,\ts\in\G_y;\ts^2=\g}\tr(\g,\s)
=(\dim\s)\i\sum_{\ts\in\G_y}\tr(\ts^2,\s)=\fra{|\G_y|}{\dim\s}.\endalign$$
Combining this with (a) we see that
$$\sum_{(x,\r)}\{(x,\r),(y,\s)\}\text{mult. of $E_{x,\r}$ in }V=1.\tag b$$

\head 4. Some applications of Theorem 0.2\endhead
\subhead 4.1\endsubhead
Let $A$ be a finite dimensional split semisimple algebra over a field $K$. Let 
$\Mod A$ be the category of $A$-modules of finite dimension over $K$. For 
$E'\in\Mod A$ let $A^{E'}$ be the sum of the simple two-sided ideals $I$ of $A$ 
such that $IE'\ne0$. For $E,E'$ in $\Mod A$ let 
$E_{E'}=A^{E'}E=\sum_{f\in\Hom_A(E',E)}f(E')$. We have the following result.

(a) {\it Let $E\in\Mod A$ and let $\cx\in A$. We have a canonical $K$-linear
isomorphism $\a:\cx E@>\si>>\Hom_A(A\cx,E)$. Moreover, $\cx E\sub E_{A\cx}$. We 
have $A\cx A\sub A^{A\cx}$.}
\nl
(Note that $A\cx$ is a left ideal of $A$ hence an object of $\Mod A$.) 
For $e\in\cx E$ we define $f_e:A\cx@>>>E$ by $f_e(a\cx)=ae,a\in A$; $f_e$ is 
well defined: if $a,a'\in A$ satisfy $a\cx=a'\cx$ then 
$ae-a'e=a\cx e_0-a'\cx e_0=0$ where $e=\cx e_0,e_0\in E$. Now $e\m f_e$ is a 
$K$-linear map $\a:\cx E@>>>\Hom_A(A\cx,E)$ which is clearly injective. We have
$\dim_K(\cx E)=\dim_K\Hom_A(A\cx,E)$. (We can assume that $A$ is a simple 
$K$-algebra and $E$ is a simple $A$-module. Thus we can assume that for some 
$K$-vector space $V$ of finite dimension we have $A=\End(V)$ and 
$E=V$ is viewed as an $A$-module in an obvious way. In this case the desired 
statement is easily verified.) It follows that $\a$ is an isomorphism.

We prove the second statement of (a). For $e\in\cx E$ we have $f_e(\cx)=e$. Since
$f_e\in\Hom_A(A\cx,E)$ we see that $e\in E_{A\cx}$, proving the second
statement of (a). Applying this to $E=A$ viewed as an object of $\Mod A$ under 
left multiplication we see that $\cx A\sub A_{A\cx}$. We now observe that 
$A_{A\cx}\sub A^{A\cx}$. Hence $\cx A\sub A^{A\cx}$ and
$A\cx A\sub AA^{A\cx}=A^{A\cx}$. This proves the third statement of (a).

\mpb

In the remainder of this section we assume that $W$ is a Weyl group and $*=1$.

\proclaim{Theorem 4.2} Let $M\in\Mod\fH$, $X\in\fH$ be as in 0.1. Let 
$E\in\Mod\fH$. We have canonically $XE\cong\Hom_\fH(M,E)$. Moreover, $XE\sub E_M$ 
and $\fH X\fH\sub\fH^M$ (notation of 4.1).
\endproclaim
We apply 4.1 with $K=\QQ(u)$, $A=\fH=\hfH$, $\cx=X$ and we use Theorem 0.2. The 
theorem follows.

\mpb

If $E$ is a simple object of $\Mod\fH$ then $\dim_K\Hom_{\fH}(M,E)$ is known from 
the work of Kottwitz \cite{\KO}; indeed, by \cite{\LV}, the specialization of our
$M$ at $u=1$ is (noncanonically) isomorphic to a $W$-module explicitly computed in
in \cite{\KO}. In particular, using the theorem we see that (a),(b) below hold.

(a) If $W$ is of type $A_n$ and $E$ is a simple $\fH$-module then $\dim_K(XE)=1$;
in particular, $E$ contains a canonical line.

(b) If $W$ is of type $B_n$ or $D_n$ and $E$ is a simple $\fH$-module then 
$\dim_K(XE)$ is a power of $2$ if $E$ is a special representation (see \cite
{\ORA}) and $XE=0$ if $E$ is a nonspecial representation.

\subhead 4.3 \endsubhead
Let $\ca=\ZZ[u,u\i]\sub K$. Let $\ch$ be the $\ca$-subalgebra of $\fH$ with basis 
$\{T_w;w\in W\}$; $\ch$ is the same as the $\ca$-algebra defined in 
\cite{\HEC, 3.2} except that $T_w,v$ of \cite{\HEC, 3.2} are the same as 
$u^{-l(w)}T_w,u$ of this paper. (When we refer to \cite{\HEC} we assume that $L=l$ 
as in \cite{\HEC, 15.1}.)

Let $J$ be the asymptotic Hecke algebra (over $\ZZ$) with basis $\{t_z;z\in W\}$ 
associated to $W$, see \cite{\HEC, \S18}. Let $\JJ=\QQ\ot J$, ${}_K J=K\ot\JJ$; 
these are split semisimple algebras.

Let $\{c_w;w\in W\}$ be the $\ca$-basis of $\ch$ as in \cite{\HEC, 5.2}. For 
$x,y,z$ in $W$ let $h_{x,y,z}\in\ca$ be as in \cite{\HEC, 13.1}. For $x,y$ in $W$ 
we write $x\si y$ if $x,y$ are in the same left cell. For $x\in W$ let 
$a(x)\in\NN$ be as in \cite{\HEC, 13.6}. Let $\cd\sub W$ be as in 
\cite{\HEC, 14.1}. The $K$-linear map $\psi:\fH@>>>{}_K\JJ$ given by
$c_x\m\sum_{d\in\cd,z\in W;d\si z\i}h_{x,d,z}t_z$ is a $K$- algebra 
isomorphism (see \cite{\HEC, 18.8}).

For any $\ce\in\Mod\JJ$ we set ${}_K\ce=K\ot\ce\in\Mod({}_K\JJ)$; let $\ce_u$ be
the $\fH$-module corresponding to ${}_K\ce$ under $\ph$. Let $\cm\in\Mod(\JJ)$ be 
such that $\cm_u\cong M$. 

From 4.2 we deduce the following result.

\proclaim{Corollary 4.4}Let $\ce\in\Mod\JJ$. We have 
$$\dim_K(\psi(X)({}_K\ce))=\dim_K\Hom_\fH(M,\ce_u).$$
Moreover, $\psi(X)({}_K\ce)\sub({}_K\ce)_{{}_K\cm}$.
\endproclaim

\subhead 4.5 \endsubhead
For $x,y,z$ in $W$ we have 

(a) $h_{x,y,z}=\g_{x,y,z\i}u^{a(z)}+\text{ lower powers of }u$ where 
$\g_{x,y,z\i}\in\NN$,
\nl
see \cite{\HEC, 13.6}. For $x\in W$ we have 
$u^{-l(w)}T_w=\sum_{y\in W}s_{y,w}c_y$ where 

(b) $s_{y,w}\in u\i\ZZ[u\i]$ for all $y\ne w$ and $s_{w,w}=1$. 

\proclaim{Proposition 4.6} Let $Z$ be a left cell of $W$ and let $a=a(h)$ for any 
$h\in Z$. Let $\x\in W$ be such that $\x\i\in Z$. We have
$$\psi(X)t_\x=\sum_{z\in Z}r_zt_zt_\x$$
where $r_z=u^a+\sum_{i<a}n_{i,z}u^i$ and $n_{i,z}\in\ZZ$ are zero for all but 
finitely many $i$.
\endproclaim
From the definitions we have
$$\align&\psi(X)t_\x=\sum_{w\in W}\psi(u^{-l(w)}T_w)t_\x
=\sum_{y,w\in W}s_{y,w}\psi(c_y)t_\x\\&=
\sum_{y,w,z\in W,d\in\cd;a(d)=a(z)}s_{y,w}h_{y,d,z}t_zt_\x.\endalign$$
By \cite{\HEC, 14.2}, in the last sum we can assume that $z\in Z$ and that 
$d\in Z$. Hence
$$\psi(X)t_\x=\sum_{y,w\in W,z\in Z,d\in\cd\cap Z}s_{y,w}h_{y,d,z}t_zt_\x.$$
Using 4.5(a),(b), we see that
$$\psi(X)t_\x=\sum_{y\in W,z\in Z,d\in\cd\cap Z}\g_{y,d,z\i}u^at_zt_\x
+\text{ lower powers of }u.$$
Using \cite{\HEC, 14.2}, we see that $\g_{y,d,z\i}$ is $1$ if $y=z$ and is $0$ 
otherwise. Thus we have 
$$\psi(X)t_\x=\sum_{z\in Z}u^at_zt_\x+\text{ lower powers of }u.$$ 
The proposition is proved.

\proclaim{Corollary 4.7} Let $Z,Z'$ be two left cells of $W$ such that 
$Z\cap Z'{}\i\ne\emp$. We have $\sum_{z\in Z\cap Z'{}\i}t_z\in\JJ^\cm$.
\endproclaim
Let $a=a(w)$ for any $w\in Z$. 
Let $d$ (resp. $d'$) be the unique element of $\cd\cap Z$ (resp. $\cd\cap Z'$).
From 4.2 we deduce (using $\psi$) that 
$\psi(X)t_d\in({}_K\JJ)^{{}_K\cm}=K\ot(\JJ^\cm)$. Using now 4.6 we deduce that 
$\sum_{z\in Z}(u^a+\sum_{i<a}n_{i,z}u^i)t_z\in K\ot(\JJ^\cm)$. It follows that
$\sum_{z\in Z}t_z\in\JJ^\cm$ hence $t_{d'}\sum_{z\in Z}t_z\in\JJ^\cm$. We now note 
that $t_{d'}\sum_{z\in Z}t_z=\sum_{z\in Z\cap Z'{}\i}t_z$. The corollary is proved.

\subhead 4.8\endsubhead
We now assume in addition that $W$ is of type $A_n,B_n$ or $D_n$. Then, by 
4.2(a),(b), the two-sided ideal $\JJ^\cm$ of $\JJ$ is the sum of the simple 
two-sided ideals corresponding to the various special representations of $W$. The 
dimension of this sum is equal to number of pairs of left cells $Z,Z'$ such that 
$Z\cap Z'{}\i$. Hence in this case, from 4.7 we deduce:

(a) {\it The elements $\sum_{z\in Z\cap Z'{}\i}t_z$ for various $Z,Z'$ as above 
form a $\QQ$-basis of $\JJ^\cm$.}
\nl
It follows that for any two-sided cell $c$ of $W$ and any left cell $Z$ contained 
in $c$,

(b) {\it the elements $\sum_{z\in Z\cap Z'{}\i}t_z$ (for various left cells $Z'$ 
contained in $c$) form a $\QQ$-basis of the unique left $\JJ$-submodule of 
$\op_{z\in Z}t_z$ isomorphic to the special representation of $\JJ$ associated to 
$c$.}

\subhead 4.9\endsubhead
For irreducible $W$ of exceptional type, the elements described in 4.8(a) do not 
span the $\QQ$-vector space $\JJ^\cm$. For example, if $W$ is of type $G_2$, that 
is, a dihedral group with generators $s_1,s_2$ such that $(s_1s_2)^6=1$, then (a) 
provides only $6$ elements while $\dim\JJ^\cm=8$. If we write $t_{12\do}$ instead 
of $t_{s_1s_2\do}$, $t_{21\do}$ instead of $t_{s_2s_1\do}$ and $t_\emp$ instead
of $t_{\text{unit element}}$, then the following $8$ elements form a $\QQ$-basis of 
$\JJ^\cm$:
$$t_\emp,t_1+t_{12121},t_{121},t_2+t_{21212},t_{212},t_{12}+t_{1212},
t_{21}+t_{2121},t_{121212}.\tag a$$
This, together with 4.8(a), suggests that for any $W$, $\JJ^\cm$ admits a 
$\QQ$-basis consisting of $\NN$-linear combinations of elements $t_z$.

\subhead 4.10\endsubhead
Let $M_\ca$ be the $\ca$-submodule of $M$ with basis $\{a_w;w\in\II_*\}$. Note that
the $\fH$-module structure on $M$ restricts to an $\ch$-module structure on $M_\ca$.
For any $\l\in\CC^*$ we regard $\CC$ as an $\ca$-module via $u\m\l$. We can then
form $M_\l=\CC\ot_\ca M_\ca$, $\ch_\l=\CC\ot_\ca\ch$ and $M_\l$ becomes a module 
over the $\CC$-algebra $\ch_\l$. Let $X_\l=1\ot X\in\ch_\l$ where $X$ is as in 0.1.
Now the assignment $a_z\m\sum_{x\in W}\tL^x_zT_z$ in 1.7 defines an $\ch$-linear
map $\mu_\ca:M_\ca@>>>\ch$ such that $\mu_\ca(a_1)=X$; by extension of scalars 
this gives rise to an $\ch_\l$-linear map $\mu_\l:M_\l@>>>\ch_\l$ such that 
$\mu_\l(a_1)=X_\l$. 
Now, if $\l\ne-1$, the $\ch_\l$-module $M_\l$ is generated by $a_1$; it follows
that in this case the image of $\mu_\l$ is the left ideal of $\ch_\l$ generated by
$X_\l$. From Theorem 0.2 it follows that there exists a finite subset $S_0$ of
$\CC^*$ such that $-1\in S_0$ and such that 

(a) for $\l\in\CC^*-S_0$, 
$\mu_\l:M_\l@>>>\ch_\l X_\l$ is an isomorphism of $\ch_\l$-modules.
\nl
(Examples in small rank suggest that one can take $S_0=\{1,-1\}$.)

\subhead 4.11\endsubhead
We now assume that $\l$ in 4.10 is such that $\l^2=q$ where $q$ is a power of
a prime number. We write $\l=\sqrt{q}$. 
Let $G$ be a split semisimple algebraic group defined over the finite field $\FF_q$
and let $G(\FF_q)$ the (finite) group of $\FF_q$-rational points of $G$. Let
$\cb$ be the flag manifold of $G$ and let $\cb(\FF_q)$ the set of 
$\FF_q$-rational points of $G$. 
Let $\cf$ be the vector space of functions $\cb(\FF_q)@>>>\CC$. For any
$B\in\cb(\FF_q)$ let $f_B\in\cf$ be the function defined by
$f_B(B')=\sqrt{q}^{l(w)}$ for any $B'\in\cb(\FF_q)$ such that $(B,B')$ are in
relative position $w\in W$.
Let $\cf'$ be the $\CC$-subspace of $\cf$ spanned by the functions $f_B$ for
various $B\in\cb(\FF_q)$. Note that $\cf$ has a natural linear action of $G(\FF_q)$
whose commuting algebra can be identified with $\ch_{\sqrt{q}}$. Then 
$\cf'$ is a $G(\FF_q)$-invariant space of $\cf$. Moreover we have
$\cf'=X_{\sqrt{q}}\cf$.
For each two-sided cell $c$ of $W$ we denote by $\cf_c$ (resp. $\cf'_c$) the sum 
of all simple $G(\FF_q)$-submodules of $\cf$ (resp. $\cf'$) which belong to $c$ in 
the classification of \cite{\ORA}. Note that $\cf_c$ is an 
$\ch_{\sqrt{q}}$-submodule of $\cf$ and that $\cf'_c=X_{\sqrt{q}}\cf_c$.
We have the following result.

\proclaim{Proposition 4.12} Assume that $\sqrt{q}\n S_0$. Let $a'=a(w_0w)$ where 
$w$ is any element of $c$ and $w_0$ is the longest element of $W$.

(a) We have $\dim(\cf'_c)=P_c(q)$ where $P_c\in\NN[t]$ ($t$ an indeterminate) is of 
the form $t^{a'}+\text{higher powers of }t$. Moreover, $P_c(1)$ 
is the number of involutions contained in $c$.

(b) We have $\dim(\cf')=P(q)$ where $P\in\NN[t]$ is such that $P(1)$ is the number
of involutions in $W$. 
\endproclaim
We prove (a). We can assume that $W$ is irreducible. The simple 
$\ch_{\sqrt{q}}$-modules which belong to $c$ can be indexed as in \cite{\ORA}
by a subset $I$ of $M(\G)$ (see 3.1) for a certain finite group $\G$ associated
to $c$; we write $\e_i$ for the simple $\ch_{\sqrt{q}}$-module indexed by $i\in I$
and $\r_i$ for the corresponding simple $G(\FF_q)$-module appearing in $\cf$.

We apply 4.1(a) with $A=\ch_{\sqrt{q}}$, $\cx=X_{\sqrt{q}}$, $E=\cf_c$. We see that 
$$\dim(\cf'_c)=\dim\Hom_{\ch_{\sqrt{q}}}(\ch_{\sqrt{q}}X_{\sqrt{q}},\cf_c).$$
Using 4.10(a) we deduce
$$\dim(\cf'_c)=\dim\Hom_{\ch_{\sqrt{q}}}(M_{\sqrt{q}},\cf_c)
=\sum_{i\in I}(\e_i:M_{\sqrt{q}})\dim\r_i\tag c$$
where $(\e_i:M_{\sqrt{q}})$ is the multiplicity of $\e_i$ in $M_{\sqrt{q}}$. As 
explained in the remarks after Theorem 0.2, the multiplicity
$(\e_i:M_{\sqrt{q}})$ can be obtained from \cite{\KO}; namely, 

if $|I|=2$, then $(\e_i:M_{\sqrt{q}})=1$ for $i\in I$;

if $|I|\ne2$, and $i=(x,\r)\in I$, then $(\e_i:M_{\sqrt{q}})$ is the multiplicity 
of $E_{x,\r}$ in $\k$ (see 3.1) or equivalently, the multiplicity of $E_{x,\r}$ in 
$V$ (see 3.2).

Thus, if $|I|=2$ we have
$$\dim(\cf'_c)=\sum_{i\in I}\dim\r_i;$$
if $|I|\ne2$ we have
$$\dim(\cf'_c)=\sum_{(x,\r)\in I}(\text{mult. of $E_{x,\r}$ in $V$})
\dim\r_{(x,\r)}.\tag d$$
Let $d(\e_i)\in\NN[t]$ be the fake degree of $\e_i$. 
If $|I|=2$ then by \cite{\ORA} we have 
$\sum_{i\in I}\dim\r_i=\sum_{i\in I}\d(\e_i)$.
If $|I|\ne2$ then by \cite{\ORA} we have
$\dim\r_i=\sum_{i'\in I}\{i,i'\}d(\e_{i'})$ where $\{i,i'\}$ is as in 3.3.
Introducing this in (d) we obtain
$$\dim(\cf'_c)=\sum_{(x,\r)\in I}(\text{mult. of $E_{x,\r}$ in $V$})
\sum_{(y,\s)\in I}\{(x,\r),(y,\s)\}d(\e_{y,\s})$$
Using now 3.3(b) we obtain
$$\dim(\cf'_c)=\sum_{(y,\s)\in I}d(\e_{y,\s}).\tag e$$
Here we have used the following two properties which are easily checked in each 
case.

$(\text{mult. of $E_{x,\r}$ in $V$})\ne0\imp(x,\r)\in I$;

If $(y,\s)\in I$ then the Frobenius-Schur indicator of $\s$ equals $1$.

We see that (e) holds both when $|I|\ne2$ and when $|I|=2$.
Now the first assertion of (a) follows immediately from (e); the second assertion
of (a) also follows from (a) using the fact that $d(\e_{y,s})|_{q=1}=
\dim(\e_{y,s})$ and that
$\sum_{(y,\s)\in I}\dim(\e_{y,\s})$ is equal to the number of involutions in $c$,
see \cite{\GE}.

Clearly, (b) is a consequence of (a). The proposition is proved.

\subhead 4.13\endsubhead
The proof of 4.12 shows that $\dim(\cf')$ is equal to the sum of the fake degrees 
$d(\e)$ of the various irreducible representations $\e$ of $\ch_{\sqrt{q}}$ (each 
one taken once).

\enddocument